\documentclass[11pt]{amsart}
\usepackage{fullpage}
\usepackage[centertags]{amsmath}
\usepackage{amsfonts}
\usepackage{amssymb}
\usepackage{amsthm}

\usepackage{pstricks}
\usepackage{epsfig}
\usepackage{pst-grad} 
\usepackage{pst-plot} 

\begin{document}

\newtheorem{teo}{Theorem}[section]
\newtheorem{coro}[teo]{Corollary}
\newtheorem{lema}[teo]{Lemma}
\newtheorem{fact}[teo]{Fact}
\newtheorem{ejem}[teo]{Example}
\newtheorem{obser}[teo]{Observation}
\newtheorem{rem}[teo]{Remark}
\newtheorem{ejer}[teo]{Exercise}
\newtheorem{propo}[teo]{Proposition}
\newtheorem{defi}[teo]{Definition}
\newtheorem{question}{Question}

\def\eqdef{\stackrel{\rm def}{=}}
\def\Proof{{\noindent {\em Proof:\ }}}
\def\fp{\hfill $\Box$}

\def\alex{\mbox{\sf AT}}
\def\baire{{\nat}^{\nat}}
\def\bairez{{\bf Z}^{\omega}}
\def\base{{\mathcal B} }
\def\binary{2^{< \omega}}
\def\cantor{2^{\nat}}
\def\cantorx{2^X}
\def\ca{{\mathcal A}}
\def\cc{{\mathcal C}}
\def\ci{{\,mathcal I}}
\def\cj{{\mathcal J}}
\def\ck{{\mathcal K}}
\def\cL{{\mathcal{L}}}
\def\coA{{\mathbb{A}}}
\def\coB{{\mathbb{B}}}
\def\coL{{\mathbb{L}}}
\def\coK{{\mathbb{K}}}
\def\cl#1{\overline{#1}}
\def\cofin{\mbox{\sf CoFIN}}
\def\ds{{\frak D}}
\def\esptop{(X,\top)}
\def\fsig{F_{\sigma}}
\def\fsd{F_{\sigma\delta}}
\def\fin{\mbox{\sf Fin}}
\def\filter{{\mathcal F}}
\def\afilter{\vec{\filter}}
\def\gfilter{\mathcal G}
\def\ged{G_\delta}
\def\hilbert{[0,1]^{\nat}}
\def\hm {O\check{c}(\mbox{$\tau$-closed, $\tau$})}
\def\ideal{{\mathcal I}}
\def\inte#1{\buildrel {\;\circ} \over #1}
\def\intseq{{\bf Z}^{<\omega}}
\def\ikl{[K,L]}
\def\iko{[K,O]}
\def\leqtau{\leq_{\tau}}
\def\k#1{{\mathcal K}(#1)}
\def\lx{<_{\scriptscriptstyle X}}
\def\ly{<_{\scriptscriptstyle Y}}
\def\rl{\Re_{l}}
\def\nat{\mathbb{N}}
\def\N{\mathbb{N}}
\def\nwd{\mbox{\sf nwd}}
\def\opcla{{\bf cl}}
\def\nin{\not\in}
\def\oc#1{O\check{c}({\fin},{#1})}
\def\oppoint{$\mbox{\bf op}$}
\def\oqpoint{$\mbox{\bf oq}$}
\def\ptres{\emptyset\times\mbox{FIN}}
\def\prodxa{\prod_\alpha \xa}
\def\power#1{{\mathcal P}(#1)}
\def\powerx{{\mathcal P}(X)}
\def\ppoint{$\mbox{p}^+$}
\def\pp{\mbox{$\mathbf  p^+$}}
\def\Pm{\mbox{$\mathbf  p^-$}}
\def\qpoint{$\mbox{q}^+$}
\def\qq{\mbox{$\mathbf{q}^+$}}
\def\Q{\mathbb{Q}}
\def\qed{\hfill \mbox{$\Box$}}
\def\R{\mathbb{R}}
\def\rpoint{$\mbox{rp}^+$}
\def\sd{\scriptstyle{\triangle}}
\def\seq{{\omega}^{<\omega}}
\def\sfan{S(\omega)}
\def\somega{S_\omega}
\def\soi{S(\ideal)}
\def\su{\subseteq}
\def\swc{\overline{S_{\Omega}}}
\def\Sq{\mbox{$\mathbf{p}^-$}}
\def\Rq{\mbox{$\mathbf{ws}$}}
\def\sw{S_{\Omega}}
\def\sx{\scriptscriptstyle X}
\def\sy{\scriptscriptstyle Y}
\def\tauf{\tau_{\scriptscriptstyle \filter}}
\def\tauaf{\tau_{\scriptscriptstyle \afilter}}
\def\tauc{\tau_{\scriptscriptstyle \cc}}
\def\tauv{\tau_{V}}
\def\ti{T(\ideal)}
\def\topy{\top_{\sy}}
\def\topo{topolog\II a }
\def\top{\tau}
\def\ult{\mathcal U}
\def\vkl{{\mathcal V}_{\scriptscriptstyle{K,L}}}
\def\veba#1#2{{\mathcal V}_{\scriptscriptstyle{#1,#2}}}
\def\iveba#1#2{[#1,#2]}
\def\ve{{\mathcal V}}
\def\ufilter{\mathcal U}
\def\w1{\omega_1}
\def\wpp{\mbox{$\mathbf{wp}^+$}}
\def\xa{X_\alpha}
\def\Z{\mathbb{Z}}

\def\sequen{\nat^{<\bf \scriptscriptstyle \omega}}
\def\casig{=^*}

\def\pie{\Pi_1^1}
\def\pieb{\boldsymbol{ \Pi}_1^1}
\def\pie03{{\mathbf{\Pi}}_3^0}
\def\sigb{{\bf \Sigma}_1^1}
\def\delb{\boldsymbol{ \Delta}_1^1}
\def\si12{\Sigma^{1}_{2}}
\def\del{\Delta_1^1}
\def\zig{\Sigma_1^1}
\def\pai2{\mathbf{\Pi}_2^0}
\def\sig2{{\bf \Sigma}^{0}_{2}}

\title{Selective separability  on spaces with an analytic topology}
\author{Javier Camargo \and Carlos Uzc\'ategui}
\address{\llap{*\,}Escuela de Matem\'aticas, Facultad de Ciencias, Universidad Industrial de
	Santander, Ciudad Universitaria, Carrera 27 Calle 9, Bucaramanga,
	Santander, A.A. 678, COLOMBIA.}
\email{jcamargo@saber.uis.edu.co}

\address{\llap{**\,}Escuela de Matem\'aticas, Facultad de Ciencias, Universidad Industrial de
	Santander, Ciudad Universitaria, Carrera 27 Calle 9, Bucaramanga,
	Santander, A.A. 678, COLOMBIA. Centro Interdisciplinario de L\'ogica y \'Algebra, Facultad de Ciencias, Universidad de Los Andes, M\'erida, VENEZUELA.}
\email{cuzcatea@saber.uis.edu.co}

\date{\today}

\thanks{The authors thank La Vicerrector\'ia de Investigaci\'on y Extensi\'on de la Universidad Industrial de Santander for the financial support for this work,  which is part  of the VIE project  \# C-2018-05}



\begin{abstract}
We study two form of  selective selective separability,  $SS$ and $SS^+$,  on   countable spaces with an analytic topology. We show several Ramsey type properties which imply $SS$. For analytic spaces $X$,  $SS^+$ is equivalent to have that the collection of dense sets is a $\ged$ subset of $2^X$, and  also equivalent to the existence of a weak base which is an $\fsig$-subset of $2^X$.   We study several examples of analytic spaces.
\end{abstract}

\maketitle

\section{Introduction}

In this paper we study some combinatorial properties of countable topological spaces. We will  focus on spaces  with a definable topology (\cite{todoruzca,Todoruzca2000}), that is to say, the topology of the space, viewed as a subset of $2^X$, has to be a definable set. Typically,  the topology will be assumed to be an analytic subset of $2^X$. The interplay between combinatorial properties of a space and the  descriptive complexity of the topology itself has  shown to be quite fruitful  \cite{CamargoUzca2018,shiva2017,todoruzca,Todoruzca2000,Todoruzca2014}.

A topological space $X$ is {\em selectively separable},  denoted $SS$, if for
any sequence $(D_n)_n$ of dense subsets of $X$ there is a finite
$a_n\subseteq D_n$ for all $n\in\nat$  such that $\bigcup_n a_n$ is dense in
$X$. This notion was introduced by Scheeper \cite{Scheeper99}  and  has received a lot of attention ever since (see for instance \cite{BarmanDow2011,BarmanDow2012,Bella2009,Bella_et_al2008,Bella2013,CamargoUzca2018,Gruenhage2011}). 
 Bella et al. \cite{Bella_et_al2008}  showed that every separable space with countable fan tightness is $SS$. On the other hand, Barman and Dow \cite{BarmanDow2011}  showed that every separable Fr\'echet space is also $SS$. In section \ref{ppoint spaces} we present several  combinatorial properties which implies  $SS$. 

Shibakov  \cite{shiva2017} showed a  stronger result when the topology of the space is analytic. He showed that any Fr\'echet countable space with an analytic topology has a countable $\pi$-base (and thus it is  $SS$). The existence of a countable $\pi$-base provides a characterization of a property  quite similar to $SS$, it is a property related to a game naturally associated to a selection principle. Let $G_1$ be the two player game defined as follows. Player I plays dense subsets of $X$ and Player II picks a point from the set played by I. So a run of the game consists of a sequence of pairs $(D_n,x_n)$ where $D_n$ is a dense set  played by I and  $x_n$ is  the response of II such that $x_n\in D_n$.  We say that II wins if $\{x_n:\;  n\in \nat\}$ is dense. Scheeper \cite{Scheeper99} showed that $X$ has a countable $\pi$-base if, and only if, II has a winning strategy for $G_1$.  For the property SS a similar game, denoted $G_{fin}$, is defined as before but  now player II picks a finite subset of the dense set played by I.  A space  $X$ has the property  $SS^+$, if II has a winning strategy for $G_{fin}$.  We show that, for  $X$  with an analytic topology, the game $G_{fin}$ is determined. We also show  that $SS^+$ is characterized by the existence of a  $\fsig$ weak $\pi$-base (that is, a $\fsig$ subset $\mathcal P$ of $2^ X$ such that  every set in $\mathcal{P}$ has non empty interior and every non empty open set contains a set from $\mathcal{P}$), which turns out to be also equivalent to having that the collection of dense subsets of $X$ is a $\ged$ subset of $2^X$.  We compare this notion of a weak base with that of  a $\sigma$-compactlike family introduced in \cite{BarmanDow2012}. Our  characterization of $SS^+$ allows to show  very easily that the product of two $SS^+$ spaces with analytic topology  is also $SS^+$. A result that holds in general as shown by Barman-Dow \cite{BarmanDow2012}. However our proof is different.  
We analyze a space constructed by Barman-Dow  \cite{BarmanDow2012} which is $SS$ and not $SS^+$ and show it has an analytic topology and has countable fan tightness. Finally, in the last section of the paper we present several examples of countable spaces.

\section{Preliminaries}
\label{preliminaries}

An {\em ideal} on a set $X$ is a collection $\ideal$ of subsets of
$X$ satisfying: (i) $A\su B$ and $B\in \ideal$, then $A\in \ideal$. (ii) If $A,B\in\ideal$, then  $A\cup B\in \ideal$. (iii) $X\nin \ideal$ and $\emptyset \in \ideal$.  We will always asume that an ideal contains all finite subsets of $X$.  If $\ideal$ is an ideal on $X$, then $\ideal^+=\{A\su X:\, A\nin \ideal\}$. 
 \fin\ denotes the ideal of finite subsets of the non negative integers $\nat$. We denote by $A^{<\omega}$ the collection of finite sequences of elements of $A$. If $s$ is a finite sequence on $A$ and $i\in A$, $|s|$ denotes its length and $ s\widehat{\;\;}i$ the sequence obtained concatenating $s$ with $i$.   For $s\in\binary$ and $\alpha\in \cantor$, let $s\prec \alpha$ if $\alpha(i)=s(i)$ for all $i<|s|$ and  $[s]=\{\alpha\in \cantor: \; s\prec \alpha\}$. The collection of all $[s]$  with $s\in\binary$ is a basis of clopen sets for $\cantor$.  Let $a\in \fin$ and $A\su \nat$, we denote by $a\sqsubseteq A$ if $a$ is an initial segment of $A$, i.e.,  $a=A\cap\{0,\cdots, n \}$, where $n=\max{a}$. For $A\su \nat$ and $m\in \nat$, we denote by $A/m$ the set $\{n\in A:\; m<n\}$ and by $A\restriction m$ the set $A\cap \{0, \cdots, m-1\}$.

Let $X$ be a topological space and $x\in X$. All spaces are assumed  to be regular and $T_1$.  A space is {\em crowded} if does not have isolated points. 
A collection $\base$ of non empty open sets is {\em a $\pi$-base}, if every non empty open set contains an element of $\base$. For every non isolated point $x$, we use the following ideal 
$$
 \ideal_x=\{A\su  X:\; x\nin \cl{A\setminus\{x\}}\}.
$$ 
The ideal of nowhere dense subsets of $X$ is denoted by $\nwd(X)$.  

Now we recall some combinatorial properties of ideals.  We put $A\su^*B$ if $A\setminus B$ is finite.
\begin{enumerate}
\item[({$\mathbf  p^+$})] $\ideal$ is \pp, if for every decreasing sequence $(A_n)_n$ of sets in $\ideal^+$, there is $A\in \ideal^+$ such that $A\su^* A_n$ for all
$n\in\nat$. Following  \cite{HMTU2017}, we say that $\ideal$   is $\Pm$, if for every decreasing sequence $(A_n)_n$ of sets in $\ideal^ +$ such that $A_n\setminus A_{n+1}\in \ideal$, there is $B\in \ideal^+$ such that $B\su^* A_n$ for all $n$.

The following notion was suggested by some results  in \cite{Farah2003,Filipowetal2008}.  By a scheme we understand a collection $(A_s)_{s\in\binary}$ such that $A_s=A_{s\widehat{\;\;}0}\;\cup\; A_{s\widehat{\;\;}1}$ and $A_{s\widehat{\;\;}0}\cap A_{s\widehat{\;\;}1}=\emptyset$ for all
$s\in\binary$.  An ideal  is \wpp,  if for every scheme $(A_s)_{s\in\binary}$ with $A_\emptyset\in \ideal^+$, there is $B\in\ideal^+$ and $\alpha\in\cantor$ such that $B\su^* A_{\alpha\restriction n}$
for all $n$. 

\item[(\qq)] $\ideal$ is \qq, if for every $A\in \ideal^+$ and every
partition $(F_n)_n$ of $A$ into finite sets, there is $S\in\ideal^+$
such that $S\su A$ and $S\cap F_n$  has at most one element for
each $n$. Such sets $S$ are called {\em selectors} for the partition. If we allow partitions with pieces in $\ideal$, we say that the ideal is {\em weakly selective} \Rq\  \cite{HMTU2017} (also called weakly Ramsey in 
\cite{SametTsaban2009}). Another natural variation is as follows: For every partition $(F_n)_n$ of a set $A\in\ideal^+$ with each piece $F_n$ in $\ideal$,  there is $S\in\ideal^+$ such that $S\su A$ and  $S\cap F_n$ is finite for all $n$. It is known that the last property  is  equivalent to \Pm\ (see  Theorem \ref{ppoint}).
\end{enumerate}

If  $\filter$ is a filter over $\nat$, recall that  $A\in \filter^+$ if $A\cap V\neq\emptyset$ for all $V\in \filter$. If $\filter^*$ is the dual ideal of $\filter$, then $\filter^+= (\filter^*)^+$. Thus  we say that a filter $\filter$ is \qq, if its dual ideal $\filter^*$ is \qq.

We say that an ideal $\ideal$ is {\em Ramsey},
usually denoted by $\ideal^+\rightarrow (\ideal^+)^2_2$, if for
every $A\in\ideal^+$ and every coloring $c:
A^{[2]}\rightarrow 2$, there is $B\su A$ in $\ideal^+$ homogeneous
for $c$, i.e. $c$ restricted to $B^{[2]}$ is constant.

A point $x$ of a topological space $X$ is called a {\em Fr\'echet
point}, if for every $A\in \ideal_x^+$  there is a sequence $(x_n)_n$ in $A$
converging to $x$. We will say
that $x$ is a \qq-{\em point}, if $\ideal_x^+$ is  \qq.   We say that space is a \qq-space, if every point is \qq. For each of the combinatorial notion about ideal introduced above,  we define analogously the  corresponding notion for points on a space.

Now we recall some other combinatorial properties of a topological space.

\begin{enumerate}

\item[(SS)]  $X$ is {\em selectively separable} \cite{Scheeper99} (see also \cite{Bella_et_al2008}), if for
any sequence $(D_n)_n$ of dense subsets of $X$ there is
$a_n\subseteq D_n$ finite for all $n\in\nat$  such that $\bigcup_n a_n$ is dense in
$X$.

\item[($\mbox{SS}^+$)]  Consider the following game $G_{fin}$.  Player I picks a dense set $D_n$; player II picks a finite $a_n \subseteq D_n$. Player II wins if  $\bigcup_n a_n$ is dense.
$X$ is $SS^+$  if player II has a winning strategy for $G_{fin}$ (\cite{BarmanDow2011}, see also \cite{Scheeper99}). Notice that if $X$ is not $SS$, then player I has an obvious winning strategy.  In Example \ref{ss-not-ss+} we show that the converse does not hold. 

\item[(RS)]  $X$ is {\em R-separable} \cite{Bella2009}, if for any
sequence $(D_n)_n$ of dense subsets of $X$ there is $x_n\in D_n$ for each $n\in\nat$ 
such that $\{x_n:\,n\in\nat\}$ is dense in $X$.

\end{enumerate}

A space  $X$ has the {\em countable fan tightness at $x$}, if for for every sequence $(A_n)_n$ with $x\in \cl{A_n}$
for all $n$, there are finite sets $K_n\su A_{n}$ for all $n$ such
that $x\in \cl{\bigcup_n K_n}$. This last property is known to be equivalent to \pp (see for instance \cite[Theorem 3.6]{SametTsaban2009}).

A subset $A$ of a Polish space is called {\em analytic}, if it is a
continuous image of a Polish space. Equivalently, if there is a
continuous function $f:\baire\rightarrow X$ with range $A$, where
$\baire$ is the space of irrationals.    For
instance, every Borel subset of a Polish space is analytic. A general reference for all descriptive set theoretic notions used in this paper is \cite{Kechris94}. We say
that a topology $\tau$ over a countable set $X$ is {\em analytic},
if $\tau$ is analytic as a subset of the cantor cube $2^X$
(identifying subsets of $X$ with characteristic functions)
\cite{shiva2017,todoruzca, Todoruzca2000,Todoruzca2014}, in this case we will say that $X$ is an {\em analytic space}. A regular countable space is analytic if, and only if, it is homeomorphic to a subspace of $C_p(\baire)$ (see \cite{todoruzca}). If there is a base $\base$ of $X$ such that $\base$ is a $\fsig$ (Borel) subset of
$\cantorx$, then we say that $X$ has a {\em $\fsig$ (Borel) base}. In general, if $X$ has a Borel base, then the topology of $X$ is analytic.

\section{\pp{},  \qq{} and SS}
\label{ppoint spaces}

In this section we show  the following implications:
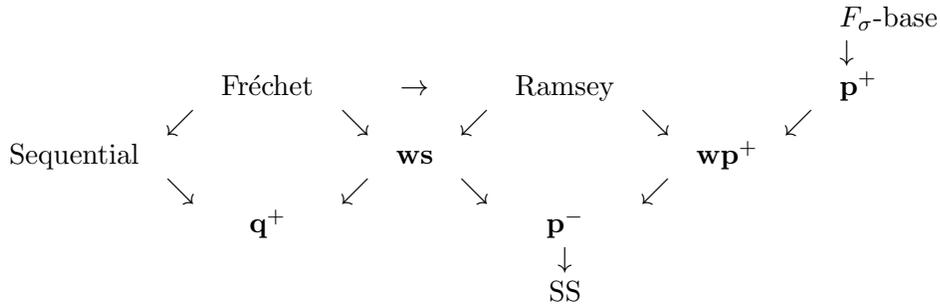
\[
\begin{array}{cccccccccclcl}
 && &                    &                          & &                  &  &                            && &  &\fsig\mbox{-base}\\
 && &                    &                          & &                  &  &                            && &  &\downarrow\\
 && &                    &\mbox{Fr\'echet}& &\rightarrow &  &\mbox{Ramsey}  && &  &\pp\\
&& &   \swarrow &&\searrow &&\swarrow& &\searrow&&\swarrow\\
&&\mbox{Sequential}&&  &                     &\mbox{\Rq}& &&  &  \mbox{\wpp} & \\
&&&\searrow & &   \swarrow &&\searrow &&  \swarrow \\
&&&&\mbox{\qq}&& & &\mbox{\Sq}\\
&&&& &&&   &\downarrow&\\
&&&& && && \mbox{SS} && 
\end{array}
\]

We start  with the ideal theoretic notions.  Some  of the  implications above are  known, however for the sake of completeness we include some of the proofs.

\begin{teo}
\label{ppoint} The following  hold for ideals on a countable set. 

\begin{itemize}

\item[(i)]  \pp\ implies \wpp.

\item[(ii)] \qq\ and \wpp\ together is equivalent to  Ramsey.

\item[(iii)] Ramsey implies  \Rq.

\item[(iv)]  \wpp\ implies \Sq.

\item[(v)] $\Pm$ is equivalent to  saying that  for every partition $(F_n)_n$ of a set $A\in \ideal^+$ with each piece $F_n$  in $\ideal$,  there is $S\in\ideal^+$ such that $S\su A$ and  $S\cap F_n$ is finite for all $n$.

\item[(vi)] \Rq\ is equivalent to \Sq\ together with \qq.

\end{itemize}
\end{teo}

\proof Let $\ideal$ be an ideal on a countable set $X$.

\begin{itemize} 
\item[(i)]  Suppose $\ideal$ is \pp.  Let  $(A_s)_{s\in\binary}$ be a
scheme with $A_\emptyset\in \ideal^+$. Let $\alpha\in \cantor$ such
that $A_{\alpha|n}\in \ideal^+$. Now we can applied \pp\ to
$(A_{\alpha|n})_n$ and finish the proof.

\item[(ii)] See \cite[Theorem 3.16]{Filipowetal2008}.

\item[(iii)] Let $(A_k)_k$ be a partition of a set $A\in \ideal^+$ such that
each $A_k\in \ideal$.  Consider the following coloring: $c(\{x,y\})=0$, if
$x, y\in A_k$ for some $k$. Let $H$ be a $c$-homogeneous set
in $\ideal^+$. Then $H$ cannot be a subset of any $A_k$, thus it has to be
1-homogeneous and thus a selector.

\item[(iv)]Suppose $\ideal$ is \wpp.  To see that $\ideal$ is \Pm, let $(A_n)_n$ be a decreasing sequence of sets in $\ideal^ +$ suth that $A_n\setminus A_{n+1}\in \ideal$.  Consider the following scheme:
\begin{enumerate}
\item[]  $B_\emptyset=A_0$,

\item[]  $B_{\langle 0\rangle}=A_0\setminus A_1$ and $B_{\langle 1\rangle}= A_1 $,

\item[] $B_{\langle 10\rangle}=A_1\setminus A_2$ and $B_{\langle 11\rangle}=A_2$,

\item[] $B_{\langle 11...10\rangle}=A_k\setminus A_{k+1}$ and $B_{\langle
11...11\rangle}=A_{k+1}$.
\end{enumerate}
and $B_s$ can be chosen arbitrarily, if $s$ is not of the form $t\widehat{\;}i$ where $t$ is a constant sequence with value 1. 

By \wpp, there is $\alpha\in\cantor$ and $C\in \ideal^+$ such that $C\su^*
B_{\alpha\restriction  m}$ for all $m$ (in fact, it is clear that the only
possibility is that   $\alpha$ is the constant sequence equal to
1). Hence $C\cap  A_n$  is finite for all $n$.

\item[(v)] This result is well known but for the sake of completeness we include a proof. Suppose $\ideal$ satisfies the second property and we show it is $\Pm$.  Let $(A_n)_n$ be a decreasing sequence of sets in $\ideal^+$ such that $A_n\setminus A_{n+1}\in\ideal$ for all $n$. Let  $F_{n}= A_n\setminus A_{n+1}$. Then $(F_n)_n$ is a partition of $A_0$ and each $F_n \in \ideal$. By hypothesis, there is  $S\in \ideal^+$ such that $S\cap F_n$ is finite. Then  clearly $S\su^*A_n$ for all $n$.  Thus $\ideal $ is $\Pm$. The other implication is shown analogously. 

\item[(vi)] This is Theorem 3.2 of \cite{HMTU2017}. We include a  proof for the sake of completeness.  We use the equivalent version of \Pm\ given in (v). Clearly \Rq\ implies \Sq\ and \qq. To check the
converse, let $A\in \ideal^+$ and $(A_n)_n$ be a partition of $A$ such
that $A_n\nin \ideal^+$ for all $n$. By \Sq\ there are finite sets
$F_n\su A_n$ such that $B=\bigcup F_n$ belongs to $\ideal^+$. From \qq, applied  to $B$ and the partition $(F_n)_n$, we get the required selector. 
\endproof
\end{itemize}

The following result is well known (see, for instance, Theorem 3.3. of \cite{HMTU2017}).

\begin{teo}
\label{fsig-p}
Every non trivial $\fsig$ ideal is \pp.
\end{teo}

\begin{coro}
\label{fsig-basis-p+} Let  $X$ be a  countable space with a $\fsig$ basis. Then every point of $X$ is  \pp\ and thus $X$ is $SS$.
\end{coro}
\proof  Suppose that $X$ has an $\fsig$ base, then it is easy to
check that $\ideal_x$ is $\fsig$ for all $x$. 

\qed

We recall that any countable subspace of $C_p(\cantor)$ admits a $\fsig$ basis (see \cite{todoruzca}).

\begin{teo} (\cite{Todoruzca2000})
\label{sec-q}
Every point of a countable sequential  space is
\qq.
\end{teo}

\begin{teo}
\label{frechet-rq} Suppose $X$ is a countable  sequential space.
\label{Frechet-w-qpoint} The following are equivalent:
\begin{itemize}
\item[(i)] $X$ is Fr\'echet.

\item[(ii)] Every point is  \Rq.

\item[(iii)] Every point is \Sq.
\end{itemize}
\end{teo}

\proof 
$(i) \Rightarrow (ii)$. Suppose $X$ is Fr\'echet.  Let $x\in X$, $A\in \ideal_x^+$  and $(F_n)_n$ be a partition of
 $A$ such that each $F_n\in\ideal_x$. We will show
that there is a selector in $\ideal_x^+$. Since $X$ is Fr\'echet, then
there is a sequence $(x_n)_n$  in $A$ converging to $x$. Since
each $F_n\in\ideal_x$, then  there are only finitely many elements
of the sequence in each $F_n$. Thus we can find a subsequence of
$(x_n)_n$ which contains at most one point of each piece $F_n$.

$(ii) \Rightarrow (iii)$. Straightforward.

$(iii) \Rightarrow (i)$. Suppose every point is  \Sq. Since $X$ is
sequential, regular and countable, it suffices to show that $X$
does not contain a copy of the Arens space. In fact, it is easy to
check that the property of being \Sq\ is hereditary and the Arens
space is not \Sq (see \S\ref{SeqSpace}).
\qed

\begin{teo} 
\label{sq-SS} Let $X$ be a countable crowded space. 
\begin{itemize}
\item[(i)] If every  point is  \Sq, then  $X$ is $SS$.

\item[(ii)] If every point  is  \Rq, then  $X$ is $RS$.

\end{itemize}
\end{teo}

\proof (i) Let $(D_n)_n$ be a sequence of dense subsets of $X$.
Let $(L_k)_k$ be a partition of $\nat$ into infinite subsets. Let
$(x_k)_k$ be an enumeration of $X$ and $(O^k_n)_{n\in L_k}$ be a maximal
family of open disjoint subsets of $X$ such that $x_k\nin
\cl{O^k_n}$ for all $n\in L_k$. Since each $D_n$ is dense, from the maximality we get that, for every $k$, 
$\bigcup_{n\in L_k} O^k_n\cap D_n$ accumulates to $x_k$. Since
each $O^k_n\cap D_n\in\ideal_{x_k}$ and $x_k$ is a \Sq-point, then
for every $k $ and $n\in L_k$ there is $E^k_n\su O^k_n\cap D_n$ finite such that  $\bigcup_{n\in
L_k}E^k_n\in\ideal_{x_k}^+$ for all $k$. Then $\bigcup_{k\in \nat}\bigcup_{n\in L_k} E^k_n$ is the
required dense set.

(ii) As before, but using \Rq\ we can assure that $|E^k_n| =1$ for all $k$ and all $n\in L_k$. 

 \qed

The following result  is an immediate corollary of Theorems  \ref{frechet-rq} and \ref{sq-SS}.

\begin{coro} (Barman-Dow \cite{BarmanDow2011})  Every countable Fr\'echet space is $RS$.
\end{coro}

Our next example shows that \pp\ does not imply $RS$. 

\begin{ejem}
\label{CL}{\em Let  $X=CL(\cantor)$ be the collection of all clopen subsets of $\cantor$ as a
subspace of $2^{\cantor}$. Then $X$ has a $\fsig$ base (see \cite{todoruzca}) and thus it is \pp and thus $SS$ (see Corollary \ref{fsig-basis-p+}). It is also known that $X$ is  not  $RS$  (see \cite{Bella2009}).   We include a proof for the sake of completeness. Let $s,t\in 2^n$, we say that $s$ and $t$ are {\em linked} if there is a sequence $u\in 2^{n-1}$ such that $s=u\widehat{\;\;}i$ and $t=u\widehat{\;\;}j$ with $i+j=1$. Each $x\in X$ is a finite union of sets of the form $[s]$ for $s\in \binary$ (see \S\ref{preliminaries}). For $k$ a positive integer,  we say that  a $x\in X$  is {\em $k$-adequated}, if   $x$ can be written as $[s_1]\cup\cdots\cup [s_m]$  with each $s_i\in 2^k$ and any pair of them are not linked. Let $A_k=\{x\in X: \; x \;\mbox{is $k$-adequated}\}$.  Notice that $\bigcup_{k\in B}A_k$ is dense in $X$ for any infinite set $B\su \nat$. Let $(B_i)_i$ be a partition of $\nat$ into infinite sets and let $D_i=\bigcup_{k\in B_i} A_k$. Then $(D_i)_i$ witnesses  that $RS$ fails. In fact,  for any selection $x_i\in D_i$, $i\in \nat$ it is easy to construct an $\alpha\in \cantor$ such that  $\alpha\nin x_i$ for all $i$.  Since $\{x\in X:\;\alpha\in x\}$ is an open set, then  $\{x_i:\; i\in \nat\}$ is not dense in $X$. 
}
\end{ejem}

\section{The property $SS^+$}
\label{SS+}

As we said in the introduction, Scheepers \cite{Scheeper99}  showed that the existence of a countable $\pi$-base is equivalent to required that player II has a winning strategy in the game $G_1$.  In this section we are going to show an analogous result for the game $G_{fin}$.

Let us say that a collection $\mathcal P$ of subsets of a space $X$ is a  {\em  weak base}  (respectively {\em  weak $\pi$-base}) if every set in $\mathcal P$ has non empty interior and for every non empty open set ·$W$ and every $x\in W$, there is $V\in \mathcal{P}$ such that $x\in int(V)\su V\su W$ (respectively,  for every non empty open set ·$W$ there is $V\in \mathcal{P}$ such that $V\su W$). Observe that if $X$ has a countable weak $\pi$-base, then it obviously has a countable $\pi$-base. The relevance of this notion is  when  we impose a complexity restriction on $\mathcal P$ (for instance, to requiere $\mathcal P$ to  be $\fsig$).  

We denote by $\ds(X)$ the collection of all dense subsets of $X$. Our first observation is the following

\begin{propo} Let $X$ be a countable  analytic space. Then 
 $\ds(X)$ is co-analytic.\qed
 \end{propo}
 
 \proof Let $\tau$ be the topology of $X$. Then $D\nin \ds(X)$ if, and only if, there is $V\in \tau$ such that $V\cap E=\emptyset$. Hence the complement of $\ds(X)$ is a projection of an analytic set and thus it is analytic. 
 \qed
 
 We present in \S\ref{complexity} an example of an analytic space $X$ such that $\ds(X)$ is not Borel. 
 
\begin{propo}
\label{fsig-pi-net} Let $Y$ be a countable set. Suppose $\ds\su 2^Y$ is closed upwards, i.e.  whenever $D\su E$ and $D\in\ds$, then $E\in \ds$.   The  following are equivalent:

\begin{enumerate}
\item[(i)] $\ds$ is $\ged$.
\item[(ii)] There are closed sets $F_n\su 2^Y$, for $n\in \nat$,  such that 
\[
2^Y\setminus \ds=\bigcup_n \{E:\;\exists V\in F_n, E\cap V=\emptyset\}.
\]
\end{enumerate}
\end{propo}

\proof It is clear that (ii) implies that the complement of $\ds$ is $F_\sigma$. Suppose now that $\ds$ is $\ged$.  Let $(K_n)_n$ be closed sets in $2^Y$  such that 
\[
\cantorx\setminus \ds=\bigcup_n K_n.
\]
Consider $F_n=\{V\su Y:\; Y\setminus V\in K_n\}$. We claim that $(F_n)_n$ satisfies  (ii). In fact,  suppose $E\nin \ds$. Then there is $n$ such that $E\in K_n$. Hence $V=Y\setminus E\in F_n$ and $V\cap E=\emptyset$. Conversely, suppose $V\in F_n$ and $E\cap V=\emptyset$. Since $Y\setminus V\nin \ds$ and $E\su Y\setminus V$, then $E\nin\ds$. 
To see that each $F_n$ is closed in $2^{Y}$ recall that the map $A\mapsto Y\setminus A$ is an homeomorphism of $2^Y$ to itself. 

\endproof

\begin{propo}
\label{fsig-weak-pi}
Let $X$ be a countable crowded space.
\begin{itemize}
\item[(i)] If $X$ has a $\fsig$ base, then $\ideal_x$ is $\fsig$ for all $x\in X$. 

\item[(ii)] If $\ideal_x$ is $\fsig$ for all $x\in X$, then $X$ has a  $\fsig$ weak base. 

\item[(iii)]  $X$ has a $\fsig$ weak $\pi$-base if, and only if,  $\ds(X)$ is $\ged$. 
\end{itemize}

\end{propo}

\proof (ii) Let $\filter_x=\{X\setminus A:\; A\in \ideal_x\}$. If  $\ideal_x$ is $\fsig$, then so is $\filter_x$. Then $\bigcup_{x\in X}\filter_x$ is a $\fsig$ weak base for $X$.

(iii) Suppose $\ds(X)$ is $\ged$ and let $(F_n)_n$ be a sequence of closed sets satisfying (ii) in  Proposition \ref{fsig-pi-net}. We claim that $\bigcup_n F_n$ is a weak $\pi$-base for $X$. First we show that if $V\in F_n$ , then $V$ has no empty interior. Since $(X\setminus V)\cap V=\emptyset$, then  by (ii), $X\setminus V\nin\ds(X)$. Thus  $V$ has non empty interior. Now let $W$ be a nonempty open set. Then $X\setminus W \nin \ds(X)$. Thus there is $n$ and $V\in F_n$ such that $(X\setminus W)\cap V=\emptyset$.  Thus $V\su W$.  
Conversely, suppose  $\bigcup_n F_n$ is a $\fsig$ weak $\pi$-base for $X$. Then it is easy to verify that  $(F_n)_n$  satisfies (ii) in Proposition \ref{fsig-pi-net}. 
\qed

As we saw in  Example \ref{CL},  $CL(\cantor)$ is a countable space with a $\fsig$ base. Since $CL(\cantor)$ does not satisfy $RS$, it cannot have  a countable $\pi$-base. 

\begin{teo}
\label{ds-gdelta}
Let $X$ be a countable crowded analytic space. Then  the game $G_{fin}$ is determined. More precisely:
\begin{itemize}
\item[(i)] If $\ds(X)$ is $\ged$, then II has a winning stratergy for $G_{fin}$.
\item[(ii)] If $\ds(X)$ is not $\ged$, then I has a winning stratergy for $G_{fin}$.
\end{itemize}
\end{teo}

\proof   We assume w.l.o.g. that $X=\nat$ and, for simplicity,  we will write  $\ds$ in place of $\ds(X)$.  For each $a\in \fin$, let $[[a]]$ be $\{A\in \cantor: a\sqsubseteq A\}$.  Notice that if $O\su\cantor$ is open and $A\in O$ with $A$ infinite, then there is $a\in \fin$ such that $A\in [[a]]\su O$. 

Suppose $\ds$ is $\ged$ and let $(O_n)_n$ be a decreasing sequence of open dense sets such that $\bigcap_n O_n ={\ds}$.  We will describe a winning strategy of player II for the game $G_{fin}$.  Let $X_0$ be the first move of player I. Since $X_0\in O_0$, player II chooses $a_0\subset X_0$ such that $[[a_0]]\su O_0$. 
Suppose $\langle X_0, a_0, \cdots X_{n-1}, a_{n-1}, X_n\rangle$ has been played such that  $X_m$ is dense for all $m\leq n$, $a_m$ is a finite subset of  $X_m$ and $[[a_0\cup\cdots\cup a_{m}]]\subseteq O_m$ for all $m<n$. Since  ·$Y=(X_n/\max{(a_0\cup\cdots\cup a_{n-1})})\cup (a_0\cup\cdots\cup a_{n-1})$ is dense in $X$, there is  a finite set $a_n\subset Y$ such that $\max(a_0\cup\cdots\cup a_{n-1})< \min a_n$ and  $[[a_0\cup\cdots\cup a_n]]\subseteq O_n$. Clearly $\bigcup_n a_n\in \bigcap_m O_m$ and hence it is dense.  This is a winning strategy for II.

Conversely, suppose $\ds$ is not $\ged$. Since $\ds$ is co-analytic,  we use  Hurewicz's theorem \cite[Theorem 21.18]{Kechris94}  to find  a countable  set $Q\su \ds$ homeomorphic to the rationals  such that  $\cl Q\cap \ds= Q$. We will describe a winning strategy for player I. Let $\{D_n:\; n\in \nat\}$ be an enumeration of $Q$. Player I starts playing $X_0= D_0$. Suppose player II plays  a finite set $a_0\subset X_0$. Since $Q$ has no isolated points, there is  $k$ such that $a_0\subset D_{k}$ and $D_{k}\neq D_0$; so  let $l_1$ be the first such $k$. Let $m_1$ be such that $D_0\restriction m_1\neq D_{l_1}\restriction m_1$ and $a_0\subset D_{l_1}\restriction m_1$. Then I plays $D_{l_1}/m_1$. 
In general,  assume the strategy has been defined up round $n$ such that if $\langle X_0, a_0, \cdots, X_{n-1}\rangle$
has been played such that player I followed this strategy, then  there are natural numbers $l_j$ and $m_j$ for $j<n$ such that, for  $i\leq n-1$,

\begin{itemize}
	\item[(i)] $X_i= D_{l_i}/m_i$,
	\item[(ii)] $a_0\cup a_1\cup\cdots\cup a_{i-1}\subset D_{l_i}\restriction m_i$,
	\item[(iii)] $D_{l_{i-1}}\restriction m_{i-1}\sqsubset D_{l_i}$, 
	\item[(iv)] $D_{l_i}\restriction m_i\neq D_j\restriction m_i$ for all $j\leq i$. 
\end{itemize}
Then II plays a finite set $a_{n-1}\subset D_{l_{n-1}}/m_{n-1}=X_{n-1}$.
For the next move, as $Q$ has no isolated points, there are  $l_n$  and $m_n$ such that $a_0\cup a_1\cup\cdots\cup a_{n-1}\subset D_{l_n}\restriction m_n$, $D_{l_{n-1}}\restriction m_{n-1}\sqsubset D_{l_n}$ and $D_{l_n}\restriction m_n\neq D_j\restriction m_n$ for all $j\leq n$.
Then  I plays $X_n=D_{l_n}/m_n$.

Finally, let  $E=\bigcup_i a_i$ and $F= \bigcup_i (D_{l_i}\restriction m_i)$. Then by (ii) $E\subseteq F$  and  by (iii) $F\in \overline{Q}$. By (iv), $F\not\in Q$ and therefore neither $F$ nor $E$ are dense.  Hence this strategy is winning for player I. 
\endproof

\begin{coro}
\label{ds-gdelta2}
Let $X$ be a countable crowded analytic space. Then $X$ is $SS^+$ if, and only if, $X$ has a $\fsig$ weak $\pi$-base. 
\end{coro}

Barman and Dow  \cite{BarmanDow2012}  introduced the notion of a compactlike family which was instrumental for their study of the property  $SS^+$. A collection $\cc$ 
of subset of a space $X$ is {\em compactlike}, if for all $D\in \ds(X)$ there is $a\su D$ finite such that $a\cap V\neq \emptyset$ for all $V\in \cc$. A space $X$ is {\em $\sigma$-compactlike}, if there is a sequence $(\cc_n)_n$ of  compactlike families of open subsets of $X$ such that the topology of $X$ is equal to $\bigcup_n \cc_n$. 

\begin{propo}
\label{CL-wbase} Let $X$ be a countable space.  Let $\cc$ be a family of subsets of $X$. Then 

(i)  $\cc$ is compactlike if, and only if, every set in $\cl\cc$ (closure in $2^X$) has non empty interior.

(ii) If $\cc$ is compactlike, then $\cl{\cc}$ is also compactlike.

(iii) If $X$ is $\sigma$-compactlike, then $X$ has a $\fsig$ weak base. 

\end{propo}

\proof (i) Suppose $\cc$ is compactlike . Let $V\in \cc$. If $V$ has empty interior, then $D=X\setminus V$ is dense. But this  contradict the definition of a compactlike family. Conversely,  suppose that every set in $\cl\cc$ has non empty interior.  Suppose $\cc$ is not compactlike and let $D$ be a dense set such that for all finite subset $a$ of $D$ there is $V\in \cc$ missing $a$. Let  $(x_n)_n$ be an enumeration of $D$ and pick $V_n\in \cc$ such that $\{x_0, \cdots, x_n\}\cap V_n=\emptyset$ for all $n$. Since $\cl\cc$ is compact, then there is $(n_k)_k$ increasing such that $V_{n_k}$ converges to some $V\in \cl\cc$.  Let $W$ be the interior of $V$ and $m$ such that $x_{m}\in W\cap D$. Thus there is $k$ such that $n_k>m$ and  $x_m \in V_{n_k}$ which is a contradiction.

(iii) Let $\cc_n$ be compactlike families of open sets such that the topology of $X$ is equal to $\bigcup_n \cc_n$. Then  by (i) and  (ii), $\bigcup_n \cl\cc_n$ is a $\fsig$ weak base for $X$. 

\qed

Barman and Dow  \cite{BarmanDow2012} showed that a countable space  satisfies $SS^+$ if, and only if, it is $\sigma$-compactlike.  Now we can show some other equivalent characterization of $SS^+$

\begin{teo}
\label{ss+charac}
Let $X$ be a countable space. The following are equivalente
\begin{itemize}
\item[(i)]  $X$ is $SS^+$.

\item[(ii)] $X$ is $\sigma$-compactlike.

\item[(iii)]  $X$ has a weak $\fsig$ base.

\item[(iv)] $X$ has a weak $\fsig$ $\pi$-base.

\item[(v)] $\ds(X)$ is $\ged$. 
\end{itemize}

\end{teo}

\proof (i)  $\Rightarrow$ (ii) was shown in  \cite{BarmanDow2012}. (ii)  $\Rightarrow$ (iii) follows from  Proposition \ref{CL-wbase}. (iii)  $\Rightarrow$ (iv) is obvious.  (iv)  $\Rightarrow$ (v) follows from Proposition  \ref{fsig-weak-pi}.  (v)   $\Rightarrow$ (i) follows from the proof of Theorem \ref{ds-gdelta}(i), just observe that in that proof the hypothesis that the topology of $X$ is analytic is not necessary. 

\qed

\medskip

Barman and Dow \cite{BarmanDow2012} have shown that the product of two countable  $SS^+$ spaces is $SS^+$. Their proof  is based in the  characterization of $SS^+$ in terms of compactlike families which makes an extensive use of ultrafilters. We present below a different proof for analytic spaces.

\begin{teo}
Let $X$ and $Y$  be countable crowded analytic spaces.  If $X$ and $Y$ are $SS^+$, then so is $X\times Y$.
\end{teo}

\proof It is clear that $X\times Y$ has an analytic topology. By  Corollary  \ref{ds-gdelta2}, it suffices to show that having a $\fsig$  weak $\pi$-base is a productive property.   Let $(F_n)_n$ and $(G_n)_n$ be closed subsets of $2^X$ and $2^Y$, respectively, such that $\bigcup_n F_n$ and $\bigcup_n G_n$ are weak $\pi$-bases for $X$ and $Y$ respectively.  Let 
$L_n^m=\{V\times W:\; V\in F_n\; \mbox{and}\; W\in G_m\}$. Then each $L_n^m$ is  a closed subset of $2^{X\times Y}$ as the map $(A,B)\in \cantorx\times 2^Y\mapsto A\times B\in 2^{X\times Y}$ is continuous. 
To finish the argument we observe that $\bigcup_{n,m} L^m_n$ is clearly a weak $\pi$-base for $X\times Y$. 
\endproof

Barman and Dow \cite[Theorem 2.12]{BarmanDow2012} constructed a countable space  $X$ which is $SS$ and not $SS^+$. We recall its definition in order to verify it has an analytic topology. We also show that it is not $SS^+$ using Theorem \ref{ds-gdelta}. They constructed  that space to show that $SS^+$ is not preserved by unions, in contrast to what happen with  $SS$  which is preserved by unions \cite{Gruenhage2011}.  In fact, there are two disjoint subspaces $A$ and $B$ such that $X=A\cup B$, $A$ and $B$ are $SS^+$, however $X$ is not $SS^+$ but it is (hereditarely) $SS$.  Moreover, $A$ and $B$ both have a $\fsig$ base.

\begin{ejem}
\label{ss-not-ss+}
{\em  We use $\alpha$, $\beta$,  etc. to denote elements of the Cantor set $\cantor$. For each $\alpha\in \cantor$, let $\alpha^*\in \cantor$ be the flipping of $\alpha$ in the first value, i.e. $\alpha^*(0)=1-\alpha(0)$ and $\alpha^* (n)=\alpha(n)$ for all $n\geq 1$.  Consider the following subset of $2^{\cantor}$
\[
Z=\{z\in 2^{\cantor}: z(\alpha)\cdot z(\alpha^*)=0\; \mbox{for all $\alpha\in \cantor$}\}.
\]
We need to modify the collection of clopen subsets of $\cantor$. Let $Q$ be the collection of eventually zero sequences. Then $Q$ is a copy of the rationals inside $\cantor$. For each $z\subseteq\cantor$, let 
\[
z_*=z\setminus Q\;\cup\;\{\alpha^*:\alpha\in z\cap Q\}.
\]
Consider the map $f$ from $2^{\cantor}$ to itself given by $f(z)=z_*$. Notice that $(z_*)_*=z$ for all $z\subseteq\cantor$. Clearly $f$ is an homeomorphism. Moreover $Z=f[Z]$.

Let $A=CL(\cantor)\cap Z$ ,   $B=\{z_*:\; z\in A\}$ and $X=A\cup B$.  Notice that  $A$ and $B$ are homeomorphic and $A\cap B=\emptyset$.   Since $CL(\cantor)$ has a $\fsig$ base  (see \cite{{todoruzca}}), then both $A$ and $B$ has a $\fsig$ base. We will show that $X$  has analytic topology, is \pp\ and  not $SS^+$.  To see that the topology of $X$ is analytic, we just need to check that the standard subbase of $X$ is analytic. In fact, for each $\alpha\in \cantor$ consider the subbasic open sets of $X$
\[
[\alpha; 1]=\{x\in X:\; \alpha\in x\}\;\;\;\mbox{and}\;\;\;  [\alpha;0]=\{x\in X:\; \alpha\nin x\}.
\]
We claim that the functions from $\cantor$ to $2^X$, $\alpha\mapsto [\alpha;i]$  are Borel measurable for $i=0,1$.  Let us check for $i=1$. Fix $x_0\in X$, it suffices to verify that $N=
\{\alpha\in \cantor:\; x_0\in [\alpha;1]\}$ is Borel. Obviously $N=x_0$. Notice that  $x_0$ is either a clopen subset of $\cantor$  or it is of the form $z_*$ for some $z\in A$ which is obviously a Borel set. Hence the collection of all subbasic open sets of the form $[\alpha;0]$ or $[\alpha;1]$ is the Borel image of the Cantor space and thus is a analytic subbase for $X$.  Therefore the topology of $X$  is also analytic (see \cite{todoruzca}). Since both $A$ and $B$ are \pp, then it is easy to verify that $X$ is also \pp.

Now we will show that $X$ is not $SS^+$ by showing that $\ds(X)$ is not $\ged$. For that end we will find a copy of the rationals $Q_1\su 2^X$ such that $\cl{Q_1}\cap \ds(X)=Q_1$. Consider
\[
Q_1=\{[q;0]\cap A\,\cup \,[q^*;0]\cap B:\; q\in Q\}.
\]
We claim that $ \cl{Q_1}=Q_1\cup \{[\alpha, 0]:\; \alpha\in \cantor\setminus Q\}$. 
We state several facts in order to verify  it.  
For notational simplicity, let 
$$
D_\beta=[\beta;0]\cap A\,\cup \,[\beta^*;0]\cap B
$$
for each $\beta\in\cantor$. 
\begin{itemize}
\item[(1)] {The map $\alpha\in \cantor\mapsto [\alpha;0]\cap A$ is continuous}:  This follows inmediately from the fact that every element of $A$ is a clopen subset of $\cantor$.

\item[(2)] {If $q_n\in Q$ and $q_n\rightarrow \alpha$ with $\alpha\nin Q$, then $[q_n;0]\cap B\rightarrow [\alpha^*;0]\cap B$}: Let $z_*\in [\alpha^*;0]\cap B$ with $z\in A$, in particular, $z$ is clopen.  Since  $\alpha^*\nin z_*$ and $\alpha\nin Q$, then $\alpha^*\nin z$. Clearly $q_n^*\rightarrow \alpha^*$, thus $q_n^*\nin z$ eventually. Therefore $q_n\nin z_*$ eventually.  That is to say $z_*\in [q_n;0]\cap B$ eventually.  Analogously, one shows that if $z_*\nin [\alpha^*;0]\cap B$, then $z_*\nin [q_n;0]\cap A$ eventually.

\item[(3)] {If $q_n\in Q$ and $q_n\rightarrow \alpha$ with $\alpha\in Q$, then $[q_n;0]\cap B\rightarrow [\alpha;0]\cap B$}:  This is shown analogously as before.

\item[(4)]  {The map $q\in Q\mapsto D_q$ is continuous}: Let $q_n\in Q$ converging to $q\in Q$. Then $D_{q_n}= [q_n;0]\cap A\,\cup \,[q_n^*;0]\cap B$ converges to $[q;0]\cap A\,\cup \,[q^*;0]\cap B$ by (1), (3) and the fact that $q_n^*\rightarrow q^*$.

\item[(5)]  Suppose $D_{q_n}\rightarrow D$ for some sequence $(q_n)_n$ in $Q$ converging to some $\alpha\in \cantor$.  
If $\alpha\in Q$, then  $D=D_\alpha$ (notice that   $D_{q_n}= [q_n;0]\cap A\,\cup \,[q_n^*;0]\cap B$ converges to $[\alpha;0]\cap A\,\cup \,[\alpha^*;0]\cap B=D_\alpha$). On the other hand, if  $\alpha\nin Q$, then  $D=[\alpha;0]$ (notice that  $D_{q_n}= [q_n;0]\cap A\,\cup \,[q_n^*;0]\cap B$ converges to $[\alpha;0]\cap A\,\cup \,[\alpha;0]\cap B=[\alpha;0]$).
%


\item[(6)] $D_\beta$ is dense for every $\beta$ and  $[\alpha;0]$ is not dense for every $\alpha$.
\end{itemize}

From (4), $Q_1$ has no isolated points and hence it is homeomorphic to the rationals.
Finally  $\cl{Q_1}\cap \ds(X)=Q_1$ by  (6). 

This example also shows  that having a winning strategy for player I in $G_{fin}$ does not imply the failure of $SS$. As this space is analytic, the game $G_{fin}$ is determined and, in fact, player I has a winning strategy because the space is not $SS^+$ (by Theorem \ref{ds-gdelta}). However, $X$ is $SS$. On the other hand, notice that in general when a space is not $SS$,  player I has an obvious winning strategy. 

}
\end{ejem}

\section{Some examples}

Now we present examples to illustrate that some of implications  shown in the previous section are strict.  Sequentiality,  $SS$ and $SS^+$  are the only notions we are considering which are strictly topological (i.e. they are no reduced to a combinatorial property of the ideal $\ideal_x$).   We will present examples showing that $SS^+$, \Pm\ and \qq\ are independent, i.e. all boolean combination of them are realizable by analytic spaces. 

 It is well known that filters (or dually, ideals) are viewed as spaces with only one non isolated point. We recall this basic construction. Suppose $Z=\nat\cup \{\infty\}$ is a space such that  $\infty $ is the only accumulation point.  Then  $\mathcal{F}_\infty=\{A\su \nat:\; \infty \in \mathrm{int}_Z(A\cup \{\infty\} )\}$ is  the neighborhood filter of $\infty$. 
Conversely, given a filter $\mathcal{F}$ over $\nat$,  we define  a topology on $\nat\cup \{\infty\}$  by declaring that  each $n\in \nat$ is isolated and $\mathcal{F}$ is the neighborhood filter of $\infty$.  We denote this space by $Z(\mathcal{F})$.    It is clear that the combinatorial properties of $\filter$ and $Z(\filter)$ are the same. Since in $Z(\filter)$ all points except $\infty$ are isolated,  we use other methods to associate  a crowded space to a filter. We  give two different such constructions. The first one is defined on $\sequen$ and is a generalization of Arkhanglel'ski{\~{\i}}-Franklin space. The second one uses inverse limits.

We need some general facts about analytic ideals and  topologies.

\begin{teo}
	\label{jalali}
	(Jalali-Naini, Talagrand \cite[Theorem 1, p. 32]{Todor97})
	Let $\ideal$ be an ideal   over a countable set  $X$ containing all finite subsets of $X$.  Suppose $\ideal$  has the Baire property as a subset of $2^X$. Then there is a partition $(K_n)_n$ of $X$  into finite sets such that $\bigcup_{i\in A}  K_i \not\in\ideal$ for all  infinite $A\subseteq X$.  
\end{teo}

\begin{propo}
\label{talagrand-general}
Let $X$ be a countable analytic crowded space and $A\su X$. For each partition $(F_n)_n$ of $A$ into finite sets there is a coarser partition $(E_n)_n$ such that $\bigcup_{n\in M} E_n $ is dense in $A$ for every infinite $M\su \nat$.
\end{propo}

\proof For each $x\in \cl A$, the ideal $\ideal_x$ (restricted to $A$) is analytic and thus it has the Baire property. By Theorem \ref{jalali},  there is a partition $(G^x_n)_n$ of $A\setminus\{x\}$ into finite sets such that $x\in \cl{\bigcup_{n\in M} G^x_n}$ for every infinite $M\su \nat$. Fix a point $x\in \cl A$. Now we define by recursion a new partition $(E^x_n)_n$ with the property for any infinite $M\su
\nat$, $\bigcup_{i\in M}E_i$ accumulates to $x$. 

We will omit the superscript $x$. Let $a_0$ be a finite set such that
$G_0\su \bigcup_{n\in a_0} F_n$. Pick $n_1$  such that $G_{n_1}\cap (\bigcup_{n\in a_0} F_n)=\emptyset$ and let $a_1$ be a finite set disjoint from $a_0$ and such that $G_{n_1}\su \bigcup_{n\in a_1} F_n$. Then we construct a sequence of disjoint finite sets $(a_i)_i$ and a sequence of integers $(n_i)_i$ such that $G_{n_i}\su \bigcup_{n\in a_i} F_n$ and $G_{n_i}\cap \bigcup_{n\in a_{i-1}} F_n=\emptyset$. Let $(k_i)_i$ be an enumeration of $\nat\setminus\bigcup_k a_k$. Let $E_i=F_{k_i}\cup \bigcup_{n\in a_i}F_n$. Then $(E_i)_i$ is a partition with the property that for any infinite $M\su
\nat$, $\bigcup_{i\in M}E_i$ accumulates to $x$. 

For each $x\in \cl A$, fix a partition $(E^x_n)_n$ (coarser than $(F_n)_n$) as before. 
The required partition $P=(E_n)_n$ can be constructed recursively so that for each $x\in \cl A$, $P$ contains infinitely many of the  blocks of the partition  $(E^x_n)_n$.
\qed

\bigskip

%
%
%
%
%
%
%
\subsection{Ideals}
We present some example of ideals satisfying some of the properties we have considered.  Probably they are known, we have included them to ilustrate some of the combinatorial properties and also because we will use them later to construct crowded spaces with similar properties.

\begin{ejem}{\em We present an example of  a \Rq\  and not Ramsey ideal.
It is well known that $\nwd(\Q)^+$ is not Ramsey (see for instance
\cite{Filipowetal2008}). 
We include the proof for the sake of
completeness.  Let $(x_n)_n$ be an enumeration of $\Q$ and define a coloring as follows:
$c: \Q^{[2]}\rightarrow \{0,1\}$, for $r<q$ in
$\Q$ put $c\{q,r\}=0$ iff $q=x_n$, $r=x_m$ and $n<m$. Then
any homogeneous set for $c$ is the range of a strictly
monotone subsequence of $(x_n)_n$ and hence is nowhere dense.

We show that $\nwd(\Q)^+$ is \Rq. Let $A\nin \nwd(\Q)$ and $(F_n)_n$ be a
partition of $A$ with each $F_n$ a nwd set. Let $(V_n)_n$ be an
enumeration of an open basis for $W=int(\cl A)$. By recursion we define a sequence of points $(x_n)_n$ and a sequence of
integers $(n_k)_k$. Pick $x_0\in V_0\cap A$ and let $n_0$ be such
that $x_0\in F_{n_0}$. Since $E_0=\bigcup_{k=0}^{n_0} F_k$ is nwd
and $V_1\cap A$ is not nwd, then pick $x_1\in (V_1\cap A)\setminus
E_0$. Let $n_1$ be such that $x_1\in F_{n_1}$. As
$E_1=\bigcup_{k=0}^{n_1} F_k$ is nwd and $V_2\cap A$ is not nwd,
then pick  $x_2\in (V_2\cap A)\setminus E_1$ and $n_2$ such that
$x_2\in F_{n_2}$ and so on. The set $S=\{x_n:\; n\in\nat\}$ is a
selector for $(F_n)_n$ which is dense in $W$. 

\qed }

\end{ejem}

\begin{ejem}
\label{fubini}{\em  The Fubini product $\ideal\times\cj$ of two ideals $\ideal$ and $\cj$ over $\nat$ is the ideal over $\nat\times\nat$ given by
\[
A\in \ideal\times\cj \Leftrightarrow \{n:\{m:(n,m)\in A\}\nin \cj\}\in \ideal.
\]
It is easy to verify that $\ideal\times \cj$ is not $\Pm$ for all ideals $\ideal$ and $\cj$. In fact, consider the sets $A_n=\{n\}\times\nat$ for $n\in \nat$. Then each $A_n\in \ideal\times \cj$, $\nat\times\nat=\bigcup_n A_n$ and for any $S\su\nat\times\nat$ such that $S\cap A_n$ is finite for all $n$, we have that $S\in \ideal\times \cj$. Thus $\ideal\times\cj$ is not $\Pm$ (see Proposition \ref{ppoint}). 

If $\cj$ is not \qq, then $\ideal\times \cj$ is not \qq\ for any ideal $\ideal$.  In fact, let $A\nin\cj$ and $(L_m)_m$ be a partition of $A$ into finite sets such that any selector for $(L_m)_m$ belongs to $\cj$. Consider the set $B=\nat\times A$. It is clear that $B\nin\ideal\times \cj$ and $B=\bigcup_{n,m} \{n\}\times L_m$. Let $S\subseteq B$ be such that $S\cap (\{n\}\times L_m)$ is finite for all $n$ and  $m$.   We claim that $S\in \ideal\times\cj$. In fact, otherwise there is $n$ such that $S'=\{m: (n,m)\in S\} \nin \cj$ and  $S'$ is a selector for $(L_m)_m$, which is a contradiction. 

We show below in Theorem \ref{qq-fubini}  that the  converse is also true for analytic ideal, i.e. $\ideal\times\ideal$ is \qq, when $\ideal$ is an analytic \qq\ ideal. 
}
\end{ejem}

\begin{ejem}
\label{no-q} {\em We recall  a well known example of a \pp\ but not \qq\ ideal. We include it for the sake of completeness. 
 Let $(K_n)_n$ be a partition of $\nat$ into finite sets such that
$|K_n|< |K_{n+1}|$ for all $n$. Let $\ideal$ be the collection of all
$A\su\nat$ such that $\sup_n |A\cap A_n|<\infty$.  Then $\ideal$ is an
ideal and clearly the partition
$(K_n)_n$ shows that $\ideal$ is not \qq. It is easy to see that $\ideal$ is an $F_\sigma$ ideal and therefore it is \pp.
\qed}
\end{ejem}

\subsection{The space $Seq(\filter)$}
\label{SeqSpace}

In this section we present a method to construct non $SS$ spaces.   We use a well-known construction of a family of topologies on $\sequen$, following  the
presentation given in \cite{Todoruzca2000}.  These topologies are defined using filters $\filter$ 
over $\nat$.   Define a topology $\tauf$ over $\sequen$
by letting a subset $U$ of $\sequen$ be open if and only if
\[
\{n\in \nat: \; s\widehat{\;}n \in U\}\in \filter
\;\;\mbox{for all $s\in U$.}
\]
Let $Seq(\filter)$ denote the space $(\sequen,\tauf)$. This space  is
$T_2$, zero dimensional and has no isolated points. Notice that
when $\filter$ is the filter of co-finite sets, then  $Seq(\filter)$ is homeomorphic to the Arkhanglel'ski{\~{\i}}-Franklin space $Seq$ and thus Arens space is homeomorphic to $\nat^{\leq 2}$ as a subspace of $Seq$.  It is also clear that $\tau_\filter$ is analytic if $\filter$ is analytic (see \cite{Uzca03} for more descriptive set theoretic properties of this space). Now we characterize the closure operator of
$Seq(\filter)$. Let $A\su \sequen$ and  put
\[
c(A)=A\cup \{s\in\sequen:\; \{n\in\nat:\; s\widehat{\;}n\in A\}\in
\filter^+\}.
\]
We can iterate this operator and define  $c_\alpha(A)$ for all
$\alpha<\omega_1$ as follows: $c_{\alpha+1}(A)=c(c_\alpha(A))$ and
$c_{\lambda}(A)=\bigcup_{\alpha<\lambda}c_{\alpha}(A)$ for
$\lambda$ limit. Then
\[
\cl{A}=\bigcup_{\alpha<\omega_1}c_\alpha(A).
\]
For $s\in \cl A$,  we define $rk(s,A)$ to be the smallest $\alpha$ such that $s\in c_\alpha(A)$.

\begin{lema}
\label{rangominimo}
For each $A\su \sequen$ and every $s\in\cl A$, there is $D\su A$ such that $s\in \cl D$ and $rk(t,E)\leq rk(s,A)$ for all $E\su D$ and every $t\in\cl E$.  
\end{lema}

\proof By induction on $rk(s,A)$. If $rk(s,A)=0$, then  let $D=\{ s\}$. Suppose $\alpha=rk(s,A)>0$ and the conclusion holds for sets with rank  smaller than $\alpha$.  By the definition of the rank, we know that 
\[
H=\{n\in\nat:\; s\widehat{\;}n\in \bigcup_{\beta< \alpha}c_\beta(A) \}\in \filter^+.
\]
Let $A_n=N_{s\widehat{\;}n}\cap A$ for $n\in H$. Notice that if $n\in H$, then $rk(s\widehat{\;}n, A)=rk(s\widehat{\;}n, A_n)<\alpha$. For each $n\in H$,  by inductive hypothesis applied to $s\widehat{\;}n$ and $A_n$, there is $D_n\su A_n$ such that $s\widehat{\;}n\in \cl{D_n}$ and $rk(t, E)\leq rk(s\widehat{\;}n, A_n)$ for all $E\su D_n$ and every  $t\in \cl{E}$. 
Let $D$ be the union of the $D_n$'s. Then $D$ satisfies the conclusion.

\qed

\begin{teo}
\label{Seq-Pm}
$Seq(\filter)$ is not SS for any filter $\filter$. In particular,  $Seq(\filter)$ is not \Sq.
\label{seq-non-SS}
\end{teo}

\proof Consider $D_n=\{s\in \sequen:\, |s|\geq n\}$. Then each $D_n$ is dense
in $Seq(\filter)$ but for any $F_n\su D_n$ finite, the set $\bigcup_n F_n$
is not dense in $Seq(\filter)$. 

\qed

The previous result implies that, among the properties we are considering,  the  only ones that  $Seq(\filter)$ could possible satisfy are sequentiality and \qq. 
When this  space is sequential  have been already characterized as stated in the following theorem. Recall that a filter $\filter$ is  {\em Fr\'echet} if for all $A\in \filter^+$ there is $B\su A$ such that every infinite subset of $B$ belongs to $\filter^+$.  

\begin{teo}
(\cite{Todoruzca2000})
\label{Seq-secuencial}  $Seq(\filter)$ is sequential if, and only if, $\filter$ is Fr\'echet. 
\end{teo}

 In \cite{GuevaraUzcategui2018} were constructed  $\aleph_1$ non homeomorphic spaces of the form $Seq(\filter)$ with $\filter$ a Fr\'echet analytic filter.
 
 Now we address the question of when $Seq(\filter)$ is \qq.
 
\begin{teo}
\label{ejem-non qpoint space} Let $\filter$ be a analytic filter on $\nat$. Then  $Seq(\filter)$ is a \qq\ if, and only if,   $\filter$ is \qq.
\end{teo}

\proof Suppose that $\filter$ is not  \qq. We will show that
$Seq(\filter)$ is not \qq\ at $\emptyset$. Let  $A\in\filter^+$ and
$(F_n)_n$ be a partition of $A$ into finite sets such that every
selector  belong to $\filter^*$. Let $B=\{s\widehat{\;}n:\; n\in
A\}$. Clearly $\emptyset\in \cl B$ and the corresponding partition
$G_m=\{s\widehat{\;}n:\; n\in F_m\}$ of $B$  does not have a
selector which accumulates to $\emptyset$. 

Conversely, suppose that $\filter$ is \qq.  For each  $s\in \sequen$, let 
\[
\ideal_s=\{A\su\sequen:\; s\nin \cl A\},
\]
where the closure is in  $Seq(\filter)$. We need to show that each $\ideal_s$  is  \qq. Let $A\su \sequen$ with $s\in\cl A$ and $(F_n)_n$  be a partition of $A$ into finite sets. Let $\alpha$ be a countable ordinal such that $s\in c_\alpha(A)$. We will show, by induction on $\alpha$, that there is a selector for $(F_n)_n$ in $\ideal_s^+$. If $s\in c(A)$, then $B=\{n\in\nat:\;s\widehat{\;}n\in A\}\in\filter^+$.  As  $(F_n)_n$ induces a partition of $B$ and $\filter$ is \qq,  there is a selector  $S\in \filter^+$ for the induced partition. Then $\{s\widehat{\;}n:\;n\in S\}$ is the required  selector for $(F_n)_n$.

Suppose now that the result holds for every $A$ and every $s\in c_\beta(A)$ with $\beta<\alpha$. Let $s\in c_\alpha(A)$ and $(F_n)_n$ a partition of $A$.  
Since the topology of $Seq(\filter)$ is analytic, then by Proposition \ref{talagrand-general} we can assume that $\bigcup_{k\in M} F_k$ is dense in $A$ for every infinite $M\su \nat$. Let $B=\{n\in\nat:\;s\widehat{\;}n\in \bigcup_{\gamma<\alpha}c_\gamma(A)\}$, then $B\in\filter^+$.  Notice that
$s\widehat{\;}n\in \cl{A\cap N_{s\widehat{\;}n}}$ and $s\widehat{\;}n\in \bigcup_{\gamma<\alpha} c_\gamma(A\cap N_{s\widehat{\;}n})$ for every $n\in B$. By Lemma \ref{rangominimo}, by passing to a subset of $A\cap N_{s\widehat{\;}n}$,  we can assume that $s\widehat{\;}n\in \cl{A\cap N_{s\widehat{\;}n}}$ and $rk(t,E)\leq rk(s\widehat{\;}n,A\cap N_{s\widehat{\;}n})<\alpha$ for all $E\su A\cap N_{s\widehat{\;}n}$ and every $t\in\cl E$.  

 Let $(M_n)_{n\in B}$ be a partition of $\nat$ into infinite sets. Let $A_n=\bigcup_{k\in M_n} F_k$, for $n\in B$.  Since each $A_n$ is dense in $A$, then $s\widehat{\;}n\in \cl{A_n}$ for $n\in B$ and by construction, $rk(s\widehat{\;}n, A_n)<\alpha$. Let $P_n$ be the partition that $(F_m)_{m\in M_n}$ induces on  $A_n\cap N_{s\widehat{\;}n}$. By the inductive hypothesis, for each $n\in B$, there is a selector $S_n\in \ideal_{s\widehat{\;}n}^+$ for $P_n$ . Then  $\bigcup_{n\in B} S_n$ is a selector of the original partition and it belongs to $\ideal_s^+$.
\endproof

\begin{coro}
\label{qq-fubini}
Let $\ideal$ be an analytic  \qq\ ideal over $\nat$. Then $\ideal\times \ideal$ is \qq. 
\end{coro}

\proof Let $\filter$ be the dual filter of $\ideal$ and  $X=\nat^2\cup \{\emptyset\}$ as a subspace of $Seq(\filter)$.  Then $\emptyset$ is an accumulation point of $X$. It is easy to verify that 
\[
\ideal_\emptyset=\{A\su X:\emptyset\nin \cl A\}=\ideal\times\ideal.
\]
Since $Seq(\filter)$ is \qq, then so is $X$ and thus $\ideal\times\ideal$ is also \qq. 

The argument above can be easily modified to show that $\ideal\times\cj$ is \qq\ when both $\ideal$ and $\cj$ are \qq.  \qed

\begin{ejem}
\label{qnoSS}  {\em $Seq$ is \qq\ (as $\fin$ is \qq) and it is not $SS$.}
\end{ejem}

\begin{ejem}
\label{DGnotq} {\em Let $\filter$ be  a non \qq\
filter (for instance, see Example \ref{no-q}). Then by Theorems \ref{ejem-non qpoint space} and \ref{seq-non-SS},   $Seq(\filter)$  is neither \qq\ nor $SS$.}

\end{ejem}
\subsection {Inverse limits}

In this section we use inverse limits to construct, for each boolean combination of the properties \Pm\ and  \qq, a  crowded $SS^+$ analytic space satisfying them. Notice first that the rationals are \pp, \qq\ and $SS^+$  by being  first countable. 

We recall a construction and some examples  of countable  spaces  presented in \cite{CamargoUzca2018}. Let $Z$ be a space and $f:Z\to Z$ continuous. 
 As usual, we denote by  $\underleftarrow{\lim}\{Z,f\}$ the inverse limit of  the constant sequence $\{Z;f\}$ which is defined by 
\[
X_\infty=\{(x_n)_n\in Z^\N:\; x_{n}=f(x_{n+1})\; \mbox{for all $n\in\N$} \}.
\]
The projection functions $\pi_n:X_\infty\rightarrow Z$ are defined by $\pi_n(x_m)_m=x_n$. 

\begin{teo} 
\label{embedding}
Let  $\mathcal{F}$ be a filter on $\nat$ and $Z=Z(\mathcal{F})$.  Let $f:Z\rightarrow Z$ be a continuous and closed surjection such that  $f(\infty)=\infty$. Let $X_\infty=\underleftarrow{\lim}\{Z,f\}$.  Then 

\begin{itemize}
\item[(i)] Every subspace of $X_\infty$ has a countable $\pi$-base. In particular, $X_\infty$ is hereditarely $SS^+$. 
\item[(ii)] If $\filter$ is \qq, then every countable subspace of $X_\infty$ is also \qq.
\item[(iii)] If $\mathcal{F}$ is \pp, then  $X_\infty$ is \pp.

\item[(iv)] For each $n\in \nat$, pick $x_n\in \pi^{-1}_1(n)$ and let $Y=\{x_n:\; n\in \nat\}\cup \{ p\}$, where $p\in X_\infty$ is the constant sequence $\infty$. Then the map $i: Z\rightarrow X_\infty$ given by $i(n)=x_n$ and $i(\infty)=p$ is an embedding. 
\end{itemize}
\end{teo}

\proof (ii)-(iv) were proven in \cite{CamargoUzca2018}. To show (i), let $Y\su X_\infty$.  Since each $n\in \N$ is isolated in $Z$, then $\{\pi^{-1}_m(n)\cap Y:\:n, m \in \N\}$ is a $\pi$-base for $Y$.

\qed

\begin{ejem} 
\label{ejemplo1}
There exists a  crowded countable analytic space $X$  which is \qq,  $SS^+$ and not \Pm.
\end{ejem}
The space $X$ will be a subspace of an  inverse limit of the form $\underleftarrow{\lim}\{Z,f\}$.  We first define the space $Z$. Consider the following ideal on $\N\times\N$: 
\[
\fin\times\fin=\{A\su\N\times\N:\; \{n\in \N:\; \{m\in\N:\; (n,m)\in A\}\not \in \fin\}\in \fin \}.
\]
Let $Z=Z(\mathcal{F})$ where $\mathcal{F}$ is the  dual filter of $\fin\times\fin$.  Since $\fin\times \fin$ is analytic, then it is  clear that the topology of $Z$ is analytic. $Z$  is \qq, as $\fin\times\fin$ is \qq (see Corollary \ref{qq-fubini}).   $Z$ is not \Pm, as $\fin\times\fin$ is not \Pm (see Example \ref{fubini}).

Let  $f\colon Z\to Z$  be any function such that:
\begin{enumerate}
	\item[(a)] $f(\infty)=\infty$;
	\item[(b)] $f^{-1}(\{n\}\times\N)=\{n\}\times\N$, for each $n$, and $f^{-1}((n,m))$ has two points for each $(n,m)\in \N^2$.
\end{enumerate}
It is not difficult to see that $f$ is a continuous, closed and open surjection. Let $X_\infty =\underleftarrow{\lim}\{Z,f\}$.  In \cite[Example 1]{CamargoUzca2018} it was shown that  every countable subspace of  $X_\infty$ has an analytic topology,   is  \qq  and SS.  By Proposition \ref{embedding}, there is an embedding $i\colon Z\to X_\infty$ such that $|\pi_1^{-1}(n,m)\cap i(Z)|=1$, for each $(n,m)\in\N\times\N$. By (b), $X_\infty$ has no isolated points, so let $D_{n,m}\su \pi^{-1}_1(n,m)$ be a countable set without isolated points such that $|D_{n,m}\cap i(Z)|=1$, for each $(n,m)\in \N\times\N$. Let $X=\bigcup_{(n,m)\in\N\times\N}D_{(n,m)}\cup\{p\}$. Since $i(Z)\subseteq X$ and $i(Z)$ is not \Pm, then  $X$ neither is \Pm.

\bigskip

\begin{ejem}
\label{ejemplo2}
There exists a regular  crowded space $X$  with a $F_\sigma$ basis (thus \pp\ and $SS^+$) which is  not \qq.
\end{ejem}

Let $\ideal$ be the ideal given in Example \ref{no-q} associated to a partition $(K_n)_n$ of $\nat$. 
Let $Z=Z(\mathcal{I}^*)$. It is clear that the topology of $Z$ is analytic.    It is obvious that  $Z$ is not \qq.  Notice that $\ideal$ is a $F_\sigma$ subset of $\cantor$, therefore $Z$ is \pp\ (see Theorem \ref{fsig-p}).  Let  $f\colon Z\to Z$  be any surjective function such that:
\begin{enumerate}
	\item[(a)] $f(\infty)=\infty$;
	\item[(b)] $f^{-1}(K_{n})=K_{n+1}$, for each $n$, and $f^{-1}(x)$ has two points for each $x\in \N$.
\end{enumerate}
It is not difficult to see that $f$ is a continuous, closed and open map. Let $X_\infty =\underleftarrow{\lim}\{Z,f\}$. In  \cite[Example 2]{CamargoUzca2018}  it was shown  that every countable subspace of $X_\infty$ has a $F_\sigma$ basis and   thus it is SS.

By (b), $X_\infty$ has no isolated points. For each $n\in \N$, let $D_n\su \pi^{-1}_1(n)$ be  a countable set without isolated points.  Let $X=\bigcup_{n}D_n\cup\{p\}$. By   Proposition \ref{embedding}, there is an embedding from $Z$ into $X$, therefore $X$ is not \qq. 

\begin{ejem}
\label{ejemplo3}
There exists a  crowded countable  analytic space $X$  which is  $SS^+$ and  neither  \qq\ nor \Pm.
\end{ejem}

Let $\ideal$ be the ideal given in Example \ref{no-q} associated to a partition $(K_n)_n$ of $\N$. Then as we have seen $\ideal\times \ideal$ is neither \qq nor \Pm.  Let $f:\N\rightarrow \N$ be an onto map  such that $f^{-1}(K_{n})=K_{n+1}$, for each $n$, and $f^{-1}(x)$ has two points for each $x\in \N$.
Let $Z=Z(\mathcal{J}^*)$. Consider $g:Z\rightarrow Z$ given by $g(n,m)=(n,f(m))$ and $g(\infty)=\infty$.  Let  $X_\infty =\underleftarrow{\lim}\{Z,g\}$. As in the previous example, let $D_{(n,m)}\su \pi^{-1}_1((n,m))$ be  a countable set without isolated points, for each $(n,m)\in \N\times\N$.  Let $X=\bigcup_{(n,m)\in\N\times\N}D_{(n,m)}\cup\{p\}$. By   Proposition \ref{embedding}, there is an embedding from $Z$ into $X$. Therefore $X$ is the required space.

\subsection{Complexity of $\ds(X)$}
\label{complexity}

We have seen that $\ds(X)$ is co-analytic when the topology of $X$ is analytic. The aim of this section is to show that there are examples where the complexity of $\ds(X)$ is co-analytic and not Borel, i.e. it is as high as possible. 

\begin{teo} Let $\filter$ be an analytic filter. Then 
$\ds(Seq(\filter))$ is a complete co-analytic set. 
\end{teo}

\proof Suppose $\filter$ is an analytic filter, then the topology  $\tauf$ is analytic and hence   $\ds(Seq(\filter))$ is co-analytic. For simplicity,  put $\ds=\ds(Seq(\filter))$. Let $\sf Tree$ be the collection of trees on $\nat$ and {\sf WFT} the collection of well founded trees on $\nat$. It is a classical fact that {\sf Tree} is a Polish space (as a subset of $2^{\sequen}$) and {\sf WFT} is a complete co-analytic set (see \cite[Section 32B]{Kechris94}).  To see that $\ds$ is   complete co-analytic we define a  Borel reduction of $\sf WFT$ into $\ds$, that is to say, a  Borel map $G:\mbox{\sf Tree}\rightarrow 2^{\sequen}$ such that $T$ is well founded iff $G(T)\in \ds$. Consider $F:\mbox{\sf Tree}\rightarrow 2^{\sequen}$ given by
\[
F(T)=\{s\in\sequen: \; \exists t\in T\;[ \, |t|=|s|\;\&\; \forall
n<|t|, \,t(n)\leq s(n)\,\,]\}.
\]
Let  $G:\mbox{\sf Tree}\rightarrow 2^{\sequen}$  be given by $G(T)=\sequen\setminus F(T)$.
We will show that  $G$ is the required  Borel reduction.
\bigskip

\noindent {\em Claim 1:} $T$ is a well founded tree iff $F(T)$ is a well
founded tree.

\Proof Notice that $F(T)$ is also a tree and $T\su F(T)$. So it
remains to show that if $F(T)$ is not well founded, then so is $T$. Let $\alpha$
be an infinite branch of $F(T)$. Consider
\[
S=\{t\in T: \; t(n)\leq \alpha(n)\; \mbox{for all $n<\mbox{lh}(t)$}\}.
\]
Then $S$ is a finitely branching subtree of $T$  and moreover $S$
is infinite (as for all $n$  there is $t_n\in T$ such that $t_n\leq
\alpha|n$). So $S$ is not well founded.
\qed

\bigskip

\noindent
{\em Claim 2:} If $T$ is not well founded, then $\mbox{int}_{\tauf}(F(T))\neq \emptyset$.

\Proof Let $\alpha$ be an infinite branch of $F(T)$ and consider
\[
V=\{t\in\sequen:\; \alpha(i)\leq t(i)\; \mbox{for all $i<\mbox{lh}(t)$} \}.
\]
Then $V\su F(T)$ and $V\in\tauf$ (if $t\in V$ and
$|t|=m$, then $t\widehat{\;}j\in V$ for all $j\geq \alpha(m)$).
\qed

\bigskip

\noindent
{\em Claim 3:} If $S$ is a well founded tree, then $S\in \nwd(Seq(\filter))$.

\Proof Since $S$ is  a tree, then it is closed in $Seq$ and therefore it is also $\tauf$-closed.  Thus it suffices to check that
$\mbox{int}_{\tauf}(S)=\emptyset$. But this is obvious, since any $\tauf$ basic nbhd
inside $S$ will provide an infinite branch in $S$.\qed

\medskip

Finally, we check that $G$ is the requiered Borel reduction of {\sf WFT} into $\ds$. Clearly  $F$ is a Borel function and so is $G$.  Now,  if  $T\in\mbox{\sf WFT}$,   then by Claim 3, $G(T)$ contains a open dense set. On the other hand, if   $T$ is not well founded,  by Claim 2,  $G(T)$ is not dense.
\endproof

\bibliographystyle{plain}


\end{document}